\input amstex
\documentstyle{amsppt}
\magnification=\magstep1                        
\hsize6.5truein\vsize8.9truein                  
\NoRunningHeads
\loadeusm

\magnification=\magstep1                        
\hsize6.5truein\vsize8.9truein                  
\NoRunningHeads
\loadeusm

\document
\topmatter

\title
The size of exponential sums on intervals of the real line
\endtitle

\author
Tam\'as Erd\'elyi, Kaveh Khodjasteh, and Lorenza Viola
\endauthor

\address Department of Mathematics, Texas A\&M University,
College Station, Texas 77843  
\endaddress

\email terdelyi\@math.tamu.edu
\endemail

\address
Department of Physics and Astronomy, 
Dartmouth College, 6127 Wilder Laboratory, Hanover, NH 03755
\endaddress

\email kaveh.khodjasteh\@dartmouth.edu
\endemail

\address
Department of Physics and Astronomy, 
Dartmouth College, 6127 Wilder Laboratory, Hanover, NH 03755
\endaddress

\email Lorenza.Viola\@Dartmouth.edu
\endemail

\endsubjclass


\abstract
We prove that 
there is a constant $c > 0$ depending only on $M \geq 1$ and $\mu \geq 0$ such that
$$\int_y^{y+a}{|g(t)| \, dt} \geq \exp (-c/(a\delta))\,, \qquad a\delta \in (0,\pi]\,,$$
for every $g$ of the form
$$g(t) = \sum_{j=0}^n{a_j e^{i\lambda_jt}},
\qquad a_j \in {\Bbb C}, \enskip |a_j| \leq Mj^\mu\,, \enskip |a_0|=1\,,
\enskip n \in {\Bbb N} \,,$$
where the exponents $\lambda_j \in {\Bbb C}$ satisfy
$$\text{\rm Re}(\lambda_0) = 0\,, \qquad \text{\rm Re}(\lambda_j) \geq j\delta > 0\,, 
\qquad j=1,2,\ldots\,, $$
and for every subinterval $[y,y+a]$ of the real line .
Establishing inequalities of this variety is motivated by problems
in physics.
\endabstract

\endtopmatter

\head {0.} Introduction \endhead

The well known Littlewood Conjecture was solved by Konyagin [8] and
independently by McGehee, Pigno, and B. Smith [10]. Based on these Lorentz [5] worked
out a textbook proof of the conjecture.

\proclaim{Theorem 0.1}
Let $n_1, n_2, \ldots, n_N$ be distinct integers. For some absolute constant $c > 0$,
$$\int_{0}^{2\pi}{\left| \sum_{k=1}^N{e^{in_k t}} \right| \, dt} \geq c \log N\,.$$
\endproclaim

This is an obvious consequence of

\proclaim{Theorem 0.2}
Let $n_1 < n_2 < \cdots < n_N$ be integers. Let $a_1, a_2, \ldots, a_k$ be
arbitrary complex numbers. We have
$$\int_{0}^{2\pi}{\left| \sum_{k=1}^N{a_ke^{in_k t}} \right| \, dt} \geq
\frac{1}{30} \sum_{k=1}^N{\frac{|a_k|}{k}}\,.$$
\endproclaim

Pichorides, who contributed essentially to the proof of the Littlewood
conjecture, observed in [10] that the original Littlewood conjecture
(when all the coefficients are from $\{0,1\}$) would follow from a result on the 
$L_1$ norm of such polynomials on sets $E \subset \partial D$ of measure $\pi$. Namely if
$$\int_E  {\Big| \sum_{j=0}^n{z^{k_j}} \Big| \, |dz|} \geq c$$
for any subset $E \subset \partial D$ of measure $\pi$ and for any nonnegative
integers $k_0 < k_1 < \cdots <k_n$ with an absolute
constant $c > 0$, then the original Littlewood conjecture holds.
Here $\partial D$ denotes the unit circle of the complex plane and the measure of a set 
$E \subset \partial D$ is the linear Lebesgue measure of the set
$$\{t \in [-\pi,\pi): e^{it} \in E\}\,.$$
Konyagin [9] gives a lovely probabilistic proof showing
that this hypothesis fails. He does however conjecture the following:
for any {\it fixed} set $E \subset \partial D$ of positive
measure there exists a constant $c = c(E) > 0$ depending only on
$E$ such that
$$\int_E {\Big| \sum_{j=0}^n{z^{k_j}} \Big|} \, |dz| \geq c(E)\,,$$
for any nonnegative integers $k_0 < k_1 < \cdots <k_n$.
In other words, the sets $E_{\varepsilon} \subset \partial D$ of measure $\pi$
in his example where
$$\int_{E_{\varepsilon}} {\Big| \sum_{j=0}^n{z^{k_j}} \Big| \, |dz|} < \varepsilon $$
must vary with $\varepsilon > 0$.

In [2] we show, among other things, that Konyagin's conjecture holds on subarcs
of the unit circle $\partial D$.

In [7] S. G\"unt\"urk constructs certain types of near-optimal approximations of a class of
analytic functions in the unit disk by power series with two distinct coefficients. More
precisely, it is shown that if all the coefficients of the power series $f(z)$ are real and
lie in $[-\mu, \mu]$, where $\mu < 1$, then
there exists a power series $Q(z)$ with coefficients in $\{-1,+1\}$ such that
$|f(z) - Q(z)| \rightarrow 0$ at the rate $\exp(C|1-z|^{-1})$ as $z \rightarrow 1$ non-tangentially
inside the unit disk. G\"unt\"urk refers to P. Borwein, Erd\'elyi, and K\'os in [5] to see that
this type of decay rate is best possible. The special case $f \equiv 0$ yields a near-optimal
solution to the ``fair duel problem" of Konyagin, as it is described in the Introduction of [7].

In this paper we extend the polynomial inequalities of [2] to
exponential sums of the form
$$ g(t) =\sum_{j=0}^n a_j e^{i\lambda_j t}\,, \qquad a_j \in {\Bbb C},
\enskip |a_j| \leq M j^\mu, \enskip |a_0|=1, \enskip  n \in {\Bbb N}, $$
\noindent
where $\text{Re}(\lambda_0)=0$ and the exponents satisfy the ``minimum
growth condition''
$$\text{Re}(\lambda_j)\geq j\delta >0, \quad j=1,2, \ldots $$
\noindent
In addition to being interesting on its own, this extension is
motivated by physical applications in the context of decoherence
control in open quantum systems using dynamical decoupling methods
[6].  In the paradigmatic case of a single two-level quantum system
undergoing pure dephasing due to coupling to a quantum bosonic
environment, for instance, the residual decoherence error at a time
$t$ after the application of $n$ ideal ``spin-flip'' pulses at times
$0<t_1<t_2< \ldots < t_n <T$, is quantified by a decay factor of the
form
$$\chi_{\{t_j\}} =\int_0^\infty \Lambda(\omega) \vert
f_{\{t_j\}}(\omega) \vert^2 d\omega\,, \quad
f_{\{t_j\}}(\omega)=\sum_{j=0}^{n} (-1)^j
(e^{it_{j}\omega}-e^{it_{j+1}\omega})\,,$$
\noindent
where $\Lambda(\omega)$ is a real function whose details depend on
both the temperature and the density of modes at frequency $\omega$ in
the environment and, in addition, $t_0:=0$ and $t_{n+1}:=T$. Thus, the
decoherence error corresponds directly to the size of the exponential
sum $f_{\{t_j \}}(\omega)$.

Decoupling methods aim to design the ``filter function'' $f_{\{t_j
\}}(\omega)$ in such a way that the error $\chi_{\{t_j\}}$ is
minimized [13]. Under the assumption that the spectral density of the
environment (hence $\Lambda (\omega)$) vanishes for frequencies higher
than an ``ultraviolet cut-off frequency'' $\omega_c$, the decay factor
$\chi_{\{t_j\}}$ may be made small by requiring the exponential sum to
vanish perturbatively, i.e., to start its Taylor series at a
sufficiently high order $(\omega_c T)^m$. In particular, it has been
recently shown [12] that if the pulse timings are chosen according to
$t_j= T \sin^2(j \pi/(2n+2))$, cancellation of $\chi_{\{t_j\}}$ is
achieved to order $m=n$ by using $n$ pulses, the so-called Uhrig
decoupling. Physically, however, a minimum growth condition is always
imposed by the fact that the separation between any two consecutive
pulses cannot be made arbitrarily small due to finite timing
resources, thus $t_{j+1}-t_j > \tau >0$ for all $j$.  As shown in [6],
the results established here may then be used to obtain a {\it
non-perturbative} lower bound on $\chi_{\{t_j\}}$, determined solely
in terms the parameter $\omega_c\tau$.  As an additional implication
of our analysis, we find that Uhrig decoupling arises naturally as a
consequence of representing certain polynomials of degree at most
$2n+1$ in terms of Lagrange interpolation at the extreme points of the
Chebyshev polynomials $U_{2n+1}$.

\head {1.} Notation \endhead

For $N > 0$ and $\mu \geq 0$, let ${\Cal S}_N^{\mu}$ denote the collection of all
analytic functions $f$ on
the open unit half-disk $D^+ := \{z \in {\Bbb C}: |z| < 1, \enskip \text{Re}(z) > 0\}$ 
that satisfy
$$|f(z)| \leq \frac {N}{\left(1-|z|\right)^\mu}\,, \qquad z \in D^+\,.$$
In this note the value of $\mu$ will always assumed to be a nonnegative integer.
We define the following subsets of ${\Cal S}_1^{1}$.
Let
$${\Cal F}_n := \left\{f: f(x) = \sum_{j=0}^n {a_jx^j}\,, \enskip
a_j \in \{-1,0,1\}\right\}$$
and denote the set of all polynomials with coefficients from
the set $\{-1,0,1\}$ by
$${\Cal F} := \bigcup_{n=0}^{\infty}{\Cal{F}_n}\,.$$
More generally we define the following classes of M\"untz polynomials.
For $M > 0$, $\mu \geq 0$, and a sequence $\Lambda := (\lambda_j)_{j=0}^{\infty}$ of complex
numbers let
$$\Cal{K}_{M}^{\mu}(\Lambda) := \left\{f: f(x) = \sum_{j=0}^n {a_jx^{\lambda_j}}\,, \enskip
a_j \in {\Bbb C}\,, \enskip |a_j| \leq Mj^{\mu}\,, \enskip |a_0|=1\,, \enskip n \in {\Bbb N} \right\}\,.$$
Here we define the analytic function $z^{\lambda_j} := \exp(\lambda_j \log z)$ by taking 
the principal analytic branch of $\log z$ in ${\Bbb C} \setminus (-\infty,0]$.  

\head {2.} New Results \endhead

\proclaim{Theorem 2.1}
There is a constant $c > 0$ depending only on $M \geq 1$ and $\mu \geq 0$ such that
$$\int_y^{y+a}{|g(t)| \, dt} \geq \exp (-c/(a\delta))\,, \qquad a\delta \in (0,\pi]\,,$$
for every $g$ of the form
$$g(t) = \sum_{j=0}^n{a_j e^{i\lambda_jt}},
\qquad a_j \in {\Bbb C}, \enskip |a_j| \leq Mj^\mu\,, \enskip |a_0|=1\,,
\enskip n \in {\Bbb N} \,,$$
where the exponents $\lambda_j \in {\Bbb C}$ satisfy 
$$\text{\rm Re}(\lambda_0) = 0\,, \qquad \text{\rm Re}(\lambda_j) \geq j\delta > 0\,, \qquad j=1,2,\ldots\,,$$
and for every subinterval $[y,y+a]$ of the real line.
\endproclaim

Using the substitution $u = t/\delta-y-a/2$ we need to prove Theorem 2.1 only
in the case when $[y,y+a] = [-a/2,a/2]$ and $\delta = 1$. Hence, using the substitution 
$z=e^{it}$, we need to prove only the following result.

\proclaim{Theorem 2.2}
There is a constant $c > 0$ depending only on $M \geq 1$ and $\mu \geq
0$ such that
$$\int_A{|f(z)| \, |dz|} \geq \exp (-c/a)\,, \qquad a \in (0,\pi]\,,$$
for every $f$ of the form
$$f(z) = \sum_{j=0}^n{a_j z^{\lambda_j}}, \qquad a_j \in {\Bbb C},
\enskip |a_j| \leq Mj^\mu\,,
\enskip |a_0|=1\,, \enskip n \in {\Bbb N} \,,$$
where the exponents $\lambda_j \in {\Bbb C}$ satisfy
$$\text{\rm Re}(\lambda_0) = 0\,, \qquad \text{\rm Re}(\lambda_j)
\geq j > 0\,, \qquad j=1,2,\ldots\,,$$
and for every subarc $A := \{e^{it}: t \in [-a/2,a/2]\}$ of the unit
circle.
\endproclaim

\medskip

\noindent {\bf Remark 2.3.} The growth condition $\text{Re}(\lambda_j) \geq j\delta$ in Theorems 2.1 and
2.2 cannot be dropped. This can be easily seen by studying the two term exponential sums
$$g_{\lambda}(t) := 1 - \exp(i\lambda t), \qquad \lambda > 0\,.$$
Clearly,
$$\lim_{\lambda \rightarrow 0+}{\max_{[-a,a]}{|g_{\lambda}(t)|}} = 
\lim_{\lambda \rightarrow 0+}{2\sin(\lambda a)} = 0\,.$$

\medskip

\noindent {\bf Remark 2.4.} The lower bounds in Theorems 2.1 and 2.2 cannot be essentially 
improved even if we assume that the exponents 
$$0 = \lambda_0 < \lambda_1 < \lambda_2 < \cdots$$
are integers and $a_j \in \{-1,1\}$ for each $j$.
Namely, in [2] the authors proved that there are absolute constants $c_1 > 0$ and $c_2 > 0$ such that
$$\inf_{0 \neq f \in \Cal{F}}{\max_{z \in A}{|f(z)|}} \leq \exp (-c_1/a)$$
whenever $A$ is a subarc of the unit circle with arclength $\ell(A) = a \leq c_2$.

\medskip

\noindent {\bf Remark 2.5.} An explicit construction showing that Theorems 2.1 and 2.2 cannot be
essentially improved can be given by utilizing the fact that if $n$ is even then
$$f_n(z) := -z+1 + 2\,\sum_{k=1}^n{(-1)^{k} z^{d_k}}\,, \qquad
d_k := \sin^2\left( \frac{k\pi}{2n+2} \right)\,, \tag 2.1$$
has a zero at $1$ with multiplicity at least $n+1$.
Namely we prove that there is an absolute constant $c > 0$ such that
$$\inf_{g}{\max_{t \in [-a,a]}{|g(t)|}} \leq 12\,\exp(-c/a)\,,$$
for all $a \in (0,1/3]$, where the infimum is taken for all exponential sums $g$ of the form
$$g(t) = \sum_{k=0}^{n+1}{a_k e^{i\lambda_kt}}$$
with $0 = \lambda_0 < \lambda_1 < \cdots < \lambda_{n+1}$ satisfying the gap condition
$$\lambda_{k+1}-\lambda_k \geq 1\,, \qquad k=0,1,\ldots, n\,,$$
($n = 2,4,\ldots$ can be arbitrary) and with 
$$a_0 = 1, \quad a_{n+1} = -1, \qquad a_k = (-1)^k2\,, \qquad k=1,2,\ldots, n\,.$$ 
Note that in the context of dynamical decoupling theory, the exponents $d_k$ are the relative 
timings of the Uhrig protocol with $n$ pulses [12].
To see that $g_n$ has a zero at $1$
with multiplicity at least $n+1$ observe that the Lagrange interpolation formula
associated with $2n+1$ distinct points (see [3, p. 8]) reproduces any polynomial of degree
at most $2n$. In particular, choosing the nodes to be the zeros
$$\alpha_k := \cos \left( \frac{k\pi}{2n+2} \right)\,, \qquad 1 \leq k \leq 2n+1\,,$$
of the Chebyshev polynomial $U_{2n+1}$ (see [3, Section 2.1] about some of
the basic facts about Chebyshev polynomials including the symmetry of the their zeros
$$\alpha_{2n+2-k} = -\alpha_k\,, \qquad k=1,2,\ldots, n\,,$$
and $\alpha_{m+1}=0$), we deduce that
$$\split Q(x) & = \sum_{k=1}^{2n+1}{Q(\alpha_k) \,
\frac{U_{2n+1}(x)}{U_{2n+1}^{\prime}(\alpha_k)(x-\alpha_k)}} \cr
& = \sum_{k=1}^{2n+1}{Q(\alpha_k) \, \frac{(-1)^{k+1}(1-\alpha_k^2)U_{2n+1}(x)}{(2n+2)(x-\alpha_k)}} \cr
& = Q(0) \, \frac{-U_{2n+1}(x)}{(2n+2)x} + \sum_{k=1}^n{
Q(\alpha_k) \, \frac{(-1)^{k+1}(1-\alpha_k^2)2x\,U_{2n+1}(x)}{(2n+2)(x^2-\alpha_k^2)}} \,, \cr
\endsplit$$
hence
$$Q(1) = Q(0) + 2\,\sum_{k=1}^n{Q(\alpha_k) \,
\frac{(-1)^{k+1}U_{2n+1}(1)(1-\alpha_k^2)}{(2n+2)(1-\alpha_k^2)}}
= 1 + 2\sum_{k=1}^n{(-1)^{k+1}Q(\alpha_k)}$$
for every polynomial $Q$ of degree at most $2n+1$. Here we used that
$$\alpha_{2n+2-k} = -\alpha_k\,, \qquad k=1,2,\ldots, n\,,$$
and $\alpha_{n+1}=0$.
Choosing
$$Q(x) := \left( 1-x^2 \right)^m\,, \qquad m=1,2,\ldots,n\,,$$
we obtain
$$0 = 1 + 2\,\sum_{k=1}^n{(-1)^{k+1} \left( 1-\alpha_k^2 \right)^m}\,, \qquad m =1,2,\ldots,n\,.$$
That is,
$$0 = -1 + 2 \,\sum_{k=1}^n{(-1)^kd_k^m}\,, \qquad m =1,2,\ldots,n\,,$$
while the assumption that $n$ is even yields
$$0 = \sum_{k=1}^n{(-1)^kd_k^m}\,, \qquad m =0\,.$$
Thus,
$$f_n^{(m)}(1) = 0, \qquad m =0,1,\ldots,n\,, \tag 2.2$$
indeed. Now let
$$g_n(t) := f_n(e^{it})\,.$$
Then
$$g_n^{(m)}(0) = 0\,, \qquad m =0,1,\ldots,n\,,$$
and
$$|g_n^{(m)}(u)| = \left| -i^m e^{iu} + 2\sum_{k=1}^n{(-1)^k(id_k)^me^{id_ku}} \right| \leq 2n+1\,,
\qquad u \in {\Bbb R}\,, \quad m=1,2,3,\ldots \,.$$
Let $b = \displaystyle{\frac{3}{n+1}}$.
By using the well known Taylor's Remainder Theorem we have
$$\split |g_n(t)| & = \left| \frac{1}{(n+1)!} \int_0^t{g_n^{(n+1)}(u)(t-u)^n \, du} \right| \cr
& \leq \frac{1}{(n+1)!} \int_0^t{\left| g_n^{(n+1)}(u)(t-u)^n \right| \, du} \cr
& \leq \frac{1}{(n+1)!} \, \max_{u \in [-|t|,|t|]}{\,|g_n^{(n+1)}(u)|} \left| t \right|^{n+1}
\leq \, \left( \frac{et}{n+1} \right)^{n+1}(2n+1) \cr
& \leq \frac{6}{b} \left( \frac{e}{3} \right)^{3/b} \cr
\endsplit \tag 2.3$$
whenever $|t| \leq 1/b$.
Now with $\lambda_0 :=0$, $\lambda_{n+1} :=9/b^2$, and $\lambda_k := (9/b^2)d_k\,, k=1,2,\ldots,n$, let
$$G_b(t) := g_n(9t/b^2) =
-e^{i\lambda_{n+1}t} + e^{i\lambda_0t} + 2\sum_{k=1}^n{(-1)^{k+1} e^{i\lambda_kt}}\,.$$
Elementary calculus shows that the exponents $\lambda_k$ satisfy the gap condition
$$\split \lambda_{k+1} - \lambda_k & \geq \lambda_1 - \lambda_0 = (9/b^2)(d_1 - d_0) =
(9/b^2)\sin^2 \left(\frac {\pi}{2n+2} \right) \cr
& \geq (9/b^2)\left(\frac {1}{n+1} \right)^2 \geq (9/b^2)(b^2/9) \geq 1\,,  
\qquad k=0,1,\ldots,n\,, \cr \endsplit$$
and it follows from (2.3) that
$$\max_{-b/9 \leq t \leq b/9}{|G_b(t)|}  
\leq \frac{6}{b} \left( \frac{e}{3} \right)^{3/b} \leq 12\exp(-c/b)$$
with an absolute constant $c>0$. 
If $b \leq 3$ is not of the form $b = \displaystyle{\frac{3}{n+1}}$
with an even non-negative integer $n$, then we choose the largest even integer $n$ such that
$b < \beta := \displaystyle{\frac{3}{n+1}}$ and the example in the already studied case
shows that
$$\max_{-b/9 \leq t \leq b/9}{|G_{\beta}(t)|} \leq 
\max_{-\beta/9 \leq t \leq \beta/9}{|G_{\beta}(t)|} \leq 12\,\exp(-c/\beta)
\leq 12\,\exp(-c^*/b)$$ with an absolute constant $c^* := c/3 > 0$.
Choosing $b = 9a \leq 3$ we obtain our claim for all $a \in (0,1/3]$.

\medskip

\noindent {\bf Remark 2.6}
Using a slightly better lower bound for $n!$, by a straightforward modification of Remark 2.5
one can see that 
$$
\inf_g \max_{t\in[-a,a]}\vert g(t)\vert \le \exp\left({-1\over{e^2 a}}\right)\left({2\over e} + ea\right)\sqrt{{e^2+1/a}\over{2\pi}}
$$
for all $a\in(0,1/(2e^2)]$, where the infimum is taken for all exponential sums $g$ of the form
$$g(t)=\sum_{k=0}^{n+1}a_ke^{i\lambda_k t},$$
with $0=\lambda_0<\lambda_1<\cdots<\lambda_{n+1}$ satisfying the gap condition
$$\lambda_{k+1}-\lambda_k\ge1,\qquad  k=0,1,\ldots,n,$$
($n=2,4,\ldots$ can be arbitrary) and with $a_k$ as in Remark 2.5. 

To see this let $n$ be even and  consider $$
\tilde{g}_{n}(z):=1-z^{\tilde{\lambda}_{n+1}}+2\sum_{k=1}^{n}(-1)^{k}z^{\tilde{\lambda}_{k}}\,,
$$
where $\tilde{\lambda}_{k}$ are defined as $$\tilde{\lambda}_{k}:=\csc^{2}\left(\frac{\pi}{2n+2}\right) \sin^{2}\left(\frac{k \pi
}{2n+2}\right),\quad k=1,\ldots,n,$$ and satisfy the gap condition $\tilde{\lambda}_{k+1}-\tilde{\lambda}_{k}\ge1$ for $k=1,2,\ldots,n$. 
Notice that the rescaled timings $\tilde{\lambda}_{k}=d_k/d_1$, where $d_k$ are defined in Remark 2.5.
Let $a=e^{-2}/n$ and define $F_{a}(t):=\tilde{g}_n(e^{it})$. In a manner similar to Remark 2.5, we can show that
$$
\vert F_{a}(t)\vert
\le \exp\left({-1\over{e^2 a}}\right)\left({2\over e} + ea\right)\sqrt{{e^2+1/a}\over{2\pi}}\,,
$$
as long as $\vert t\vert \le a$. Notice that r.h.s of the inequality we just derived is an increasing function of $a$. If $a$ cannot be written as $e^{-2}/n$ for an even integer $n$, we may simply use the smallest even integer $n$ such that $\displaystyle a>\alpha:=e^{-2}/n$ as long as $0<a<1/(2e^2)$
and thus
$$
\max_{-a\le t\le a} \vert F_\alpha(t)\vert \le \max_{-\alpha \le t\le \alpha } \vert F_\alpha(t)\vert \le \exp\left({-1\over{e^2 a}}\right)\left({2\over e} + ea\right)\sqrt{{e^2+1/a}\over{2\pi}}\,.
$$

\medskip

\noindent {\bf Remark 2.7}
Finding a polynomial $\sum_{j=1}^n{a_jz^{\lambda_j}}$ with $a_j \in \{-1,1\}$, {\it integer} 
exponents $0 < \lambda_1 < \lambda_2 < \cdots < \lambda_{2n}$, and with a zeros at $1$
of multiplicity at least $n$ is closely related to Wright's conjecture (1934) on ideal 
solutions of the Prouhet-Tarry-Escott Problem. This seems extremely difficult to settle.
See [1, Chapter 11] about the history of this problem. However, as P. Borwein writes it 
in [1, p. 87], heuristic arguments suggest that Wright's conjecture should be false.

\medskip

\head {3.} Lemmas \endhead

To prove Theorem 2.2 we modify the proof given in [2] in the case where $\lambda_j = j$ for
each $j$. We need the lemmas.  

\proclaim{Lemma 3.1} Let $0 < a \leq \pi$ and $N \geq 1$. For every $g \in {\Cal S}_N^{\mu}$ with
$|g (\frac{1}{4Ne^{\mu}})| \geq 4^{-(\mu+1)}$ there is a value $b \in [\frac 12, \frac 34]$ 
such that $|g(b)| \geq c_2 > 0$, where $c_2$ depends only on 
$N \geq 1$ and $\mu \geq 0$.
\endproclaim

\demo{Proof of Lemma 3.1}
The proof is a standard normal family argument. Suppose the lemma is not true for some $N \geq 1$
and $\mu \geq 0$. Then there is a sequence $(g_n)$  such that $g_n \in {\Cal S}_N^{\mu}$,
$| g_n (\frac{1}{4Ne^{\mu}})| =  4^{-(\mu+1)}$ and 
$$\lim _{n \rightarrow \infty}{K_n} = 0\,, \qquad 
K_n := \max_{z \in [\frac 12, \frac 34]}{|g_n(z)|}\,, \qquad n=1,2,\ldots\,.$$
Then there is a subsequence of $(g_n)$,
without loss of generality we may assume that this is $(g_n)$ itself, that converges
to a function $g \in {\Cal S}_N^{\mu}$ locally uniformly on every compact subset of
$$H := \{z \in {\Bbb C}: |z| < 1, \enskip \text{Re}(z) > 0\}\,.$$
Now $g(z)=0$ for all $z \in [\frac 12, \frac 34]$ while 
$g( \frac{1}{4Ne^{\mu}} ) = 4^{-(\mu+1)}\,.$  This contradicts the Unicity Theorem.
\qed \enddemo

\proclaim{Lemma 3.2} Let $0 < a \leq  \pi$. $b \in [\frac 12, \frac 34]$, and $N \geq 1$. 
Let $\Gamma_{a,b}$ be the circle centered at $b$ with diameter $[2b-\cos(a/8),\cos(a/8)]$.
Let $I$ be the subarc of $\Gamma_{a,b}$ with length $\ell(I) \geq c_3a$ with midpoint 
$\cos(a/8)$ on the real line, where $c_3 \in (0,1]$.
Then there is a constant $c_4 > 0$ depending only on $N$, $\mu$, $c_2$, and $c_3$ such that
$$\max_{z \in I}{|g(z)|} \geq \exp (-c_4/a)$$ 
for every $g \in {\Cal S}_N^{\mu}$ with  $|g(b)| \geq c_2\,.$
\endproclaim

\demo{Proof of Lemma 3.2}
Let $2m \geq 4$ be the smallest even integer not less than $4\pi/(c_3a)$.
Let
$$\xi := \exp \left( \frac {\pi i}{m} \right)$$
be the first $(2m)$-th root of unity.
We define $2m$ equally spaced points on $\Gamma_{a,b}$ by
$$\eta_k := b + (\cos(a/8) - b)  \xi ^k\,, \qquad k=0,1,\ldots, 2m-1\,.$$
Then there is a constant $c_5 > 0$ depending only on $c_3$ such that
$$1 - |z| \geq c_5 (ka)^2\,, \qquad k=1,2, \ldots m-1\,,$$
whenever $z$ is on the smaller subarc of the circle $\Gamma_{a,b}$ with endpoints
$\eta_k$ and $\eta_{k+1}$ or with endpoints $\eta_{2m-k}$ and $\eta_{2m-k-1}$,
respectively. We define the function
$$h(z) := \prod_{j=0}^{2m-1}
{g(b + (\cos(a/8) - b) \xi^j (z - b))}\,.$$
If $g \in {\Cal S}_N^\mu$, then
$$\leqalignno{\max_{z \in \Gamma_{a,b}} {|h(z)|} & \leq 
\prod_{k=1}^{m-1}{\left( N\left(\frac{1}{c_5(ka)^2} \right)^{\mu} \max_{z \in I}{|g(z)|} \right)^2} \cr 
& \leq \left(\frac{1}{c_5 a}\right)^{(4m-4)\mu} \frac {N^{2m-2}}{((m-1)!)^{4\mu}}
\left( \max_{z \in I}{|g(z)|} \right)^2 \cr
& \leq \left( \frac{m}{2\pi c_3c_5} \right)^{(4m-4)\mu} \left( \frac{e}{m-1} \right)^{(4m-4)\mu} 
N^{2m-2} \left( \max_{z \in I}{|g(z)|} \right)^2 \cr 
& \leq e^{4\mu} N^{2m-2} \left( \frac{e}{2\pi c_3c_5} \right)^{(4m-4)\mu} 
\left( \max_{z \in I}{|g(z)|} \right)^2 \cr
& \leq \exp (c_6/a) \, \left( \max_{z \in I}{|g(z)|} \right)^2 \cr} $$
with a constant $c_6 > 0$ depending only on $N$, $\mu$, and $c_3$.
Now the Maximum Principle yields that
$$\left|g(b)\right|^{2m} = \left|h(b)\right| \leq \max_{z \in \Gamma_{a,b}} {|h(z)|}
\leq \exp (c_6/a)  \, \left( \max_{z \in I}{|g(z)|} \right)^2\,.$$
Since $2m \leq 2 + 4\pi/(c_3a)$ and $|g(b)| \geq c_2$,
we obtain
$$\left( \max_{z \in I}{|g(z)|} \right)^2 \geq \exp (-c_6/a) \left|g(b)\right|^{2m}
\geq \exp (-c_6/a) \left( c_2 \right)^{2m} \geq \exp (-2c_4/a) $$
with a constant $c_4 > 0$ depending only on $N$, $\mu$, $c_2$, and $c_3$. 
\qed \enddemo

\proclaim{Lemma 3.3} Let $0 < a < \pi$, $N \geq 1$, and $\mu = 1,2, \ldots$.
Let $A := \{e^{it}: t \in [-a/2,a/2]\}$. There is a constant $c_7 > 0$ depending 
only on  $N$ and $\mu$ such that
$$\int_A{|f(z)| \, |dz|} \geq \exp(-c_7/a)$$
for every $f \in {\Cal S}_N^{\mu}$ that is analytic on the arc $A$ and
satisfies $|f ( \frac{1}{4Ne^{\mu}} )| \geq \frac 12\,.$
\endproclaim

To prove Lemma 3.3 we need the following.

\proclaim {Lemma 3.4}
Let $w_1 \neq w_2 \in {\Bbb C}$ and let $z_0 := \frac 12 (w_1+w_2)$.
Assume that $J_1$ is an arc that connects
$w_1$ and $w_2$. Let $J_2$ be the arc that is the symmetric image of
$J_1$ with respect to the $z_0$.
Let $J := J_1 \cup J_2$ be positively oriented.
Suppose that $g$ is an analytic function inside and on $J$. Suppose
that the region inside $J$ contains the disk centered at $z_0$
with radius $\gamma > 0$. Let $|g(z)| \leq K$ for $z \in J_2$.
Then
$$|g(z_0)|^2 \leq (\pi \gamma)^{-1} K \int_{J_1}{|g(z)|\,|dz|}\,.$$
\endproclaim

\demo{Proof of Lemma 3.4}
Applying Cauchy's integral formula with
$$G(z) := g(z_0+(z-z_0))g(z_0-(z-z_0))$$
on $J$, we obtain
$$\leqalignno{|g(z_0)|^2 & = |G(z_0)|
= \left| \frac{1}{2\pi i} \,\int_{J} {\frac {G(z)\,dz}{z-z_0}} \right| \cr
& = \frac {2}{2\pi}\,\left| \int_{J_1} {\frac {G(z) \,dz}{z-z_0}} \right|
\leq \frac {1}{\pi} \int_{J_1} {\frac {|G(z)|\,|dz|}{|z-z_0|}} \cr
& = \frac {1}{\pi} \int_{J_1} {\frac {|g(z_0+(z-z_0))g(z_0-(z-z_0))|\,|dz|}{|z-z_0|}} \cr
& \leq (\pi \gamma)^{-1} K \int_{J_1}{|g(z)|\,|dz|}\,. \cr}$$
\qed
\enddemo

\demo{Proof of Lemma 3.3}
Without loss of generality we may assume that $\ell(A) = a \leq \pi/2$. Suppose 
$f \in {\Cal S}_N^{\mu}$ and $|f(\frac {1}{4Ne^{\mu}})| \geq \frac 12\,.$
Let the region $H_a$ be defined by
$$H_a := \left\{z = re^{i\theta}: \enskip \cos(a/4) < r < \cos(a/8)\,,
\enskip -a/4 < \theta < a/4 \right\}\,.$$
Associated with $a \in (0,1]$ and $b \in [\frac 12,\frac 34]$ (the choice of $b$
will be specified later), let 
$\Gamma_{a,b}$ be the circle as in Lemma 3.2. It is easy to see that the arc 
$I := \Gamma_{a,b} \cap H_a$ has length greater than $c_3a$ with an absolute constant 
$c_3 > 0$. Let $f \in {\Cal S}_N^{\mu}$. Let $z_0 \in I \subset H_a$ be chosen so that
$$|f(z_0)| = \max_{z \in I}{|f(z)|}\,.$$
Also, we can choose $w_1 \in A$ and
$w_2 \in A$ such that $z_0 = \frac 12 \left( w_1 + w_2 \right)$.
Let $J_1$ be the arc connecting $w_1$ and $w_2$ on the unit circle. Note that
$J_1$ is a subarc of $A$ of length at least $a/4$.
Let $J_2$ be the arc which is the symmetric image of $J_1$ with respect to
the line segment connecting $w_1$ and $w_2$.
Let
$$g(z) := 4^{-\mu}((z-w_1)(z-w_2))^\mu f(z)\,.$$
It is elementary geometry again to show that
$$|g(z)| \leq \frac {4^{-\mu}N|(z-w_1)(z-w_2)|^\mu}{\left(1-|z|\right)^\mu}
\leq \frac{4^{-\mu}N2^\mu}{\sin^\mu(a/8)}\,= \frac{2^{-\mu}N}{\sin^\mu(a/8)}\,, \qquad z \in J_2\,.$$
By Lemma 3.4 we obtain
$$|g(z_0)|^2 \leq (\pi(1 - \cos(a/8))^{-1} 
\frac {N2^{-\mu}}{\sin^\mu(a/2)} \int_{J_1}{|g(z)|\, |dz|}\,. \tag 3.1$$
Observe that $f \in {\Cal S}_N^{\mu}$ implies $g \in {\Cal S}_N^{\mu}$.
Also, since $N \geq 1$, $\mu \geq 1$, and $|f ( \frac{1}{4Ne^{\mu}} )| \geq \frac 12\,$, we have 
$$\split \left|g\left( \textstyle {\frac {1}{4Ne^\mu}} \right)\right|
\geq & 4^{-\mu} \left(1 - \textstyle{\frac{1}{4Ne^{\mu}}} \right)^{2\mu} 
\left|f\left( \textstyle {\frac {1}{4Ne^\mu}} \right)\right| 
\geq 4^{-\mu} \left(1 - \textstyle{\frac {1}{8\mu}}\right)^{2\mu} \textstyle{\frac 12} 
\geq 4^{-\mu} \textstyle{(\frac 78 )^2 \frac 12} \cr
\geq & 4^{-(\mu+1)}\,. \cr \endsplit $$
Hence, by Lemma 3.1, we can pick $b \in [\frac 12,\frac 34]$ so that $|g(b)| \geq c_2$
with an absolute constant $c_2 > 0$ depending only on $N$ and $\mu$.
Now we can deduce from Lemma 3.2 that
$$|g(z_0)| \geq \exp(c_4/a)\,. \tag 3.2$$
Combining  (3.2) with (3.1) and $J_1 \subset A$ gives
$$\leqalignno{\int_A{|f(z)| \, |dz|}\, & \geq \, \int_A{|f(z)| \, |dz|} 
\geq \int_{J_1}{|g(z)| \, |dz|} \cr 
& \geq \, \pi(1 - \cos(a/8)) \,\frac {2^{\mu}\sin^{\mu}(a/8)}{N}\, |g(z_0)|^2 \cr
& \geq \exp (-c_1/a) \cr}$$
with a constant $c_1 > 0$ depending only $N$ and $\mu$.
\qed \enddemo

\head {4.} Proof of the theorems \endhead

\demo{Proof of Theorem 2.2} Let $f \in {\Cal K}_M^{\mu}(\Lambda)$, where
where the exponents $\lambda_j \in {\Bbb C}$ satisfy
$$\text{\rm Re}(\lambda_0) = 0\,, \qquad \text{\rm Re}(\lambda_j) \geq j > 0\,, \qquad j=1,2,\ldots\,.$$
Then $f \in {\Cal S}_{M\mu!}^{\mu+1}$ and $f$ is analityc on the arc $A$.
Also, if $|z_0| \leq \textstyle {\frac {1}{4M(\mu+1)!e^{\mu+1}}}$, then
$$\leqalignno{|f(z_0)| & \geq 1 - M \, \sum_{j=1}^{\infty} 
{\lambda_j^{\mu}\left(\frac {1}{4M(\mu+1)!e^{\mu+1}} \right)^{\lambda_j}} 
\geq 1 - \frac {M}{4M(\mu+1)!}\sum_{j=1}^{\infty} {\left(\frac {j}{e^j} \right)^{\mu+1}} \cr
& \geq 1 - \frac 14 \sum_{j=1}^{\infty}{\left(\frac {j}{e^j} \right)} \geq
1 - \frac 14 \sum_{j=1}^{\infty}{\left(\frac {j}{2^j} \right)} 1 - \frac 24 \geq \frac 12\,. \cr}$$
So the assumptions of Theorem 3.3 are satisfied with $N$ replaced by
$M\mu!$, and the theorem follows from Lemma 3.3. 
\qed
\enddemo

Theorem 2.1 is an obvious consequence of Theorem 2.2.

\head
{Acknowledgments}
\endhead

Work at Dartmouth was supported by the National Science Foundation
under grant number PHY-0903727.

\enddocument